\documentclass{gtmon_a}
\pdfoutput=1
\usepackage{graphicx}

\proceedingstitle{Heegaard splittings of 3--manifolds (Haifa 2005)}
\conferencestart{10 July 2005}
\conferenceend{19 July 2005}
\conferencename{Heegaard splittings of 3--manifolds}
\conferencelocation{Haifa}

\editor{Cameron Gordon}
\givenname{Cameron}
\surname{Gordon}

\editor{Yoav}
\givenname{Yoav}
\surname{Moriah}

\title[Splittings of $S^3$ are unique]{A proof of Waldhausen's
uniqueness of splittings of $S^3$ (after Rubinstein and
Scharlemann)}
\author{Yo'av Rieck}%
\givenname{Yo'av}
\surname{Rieck}
\address{Department of mathematical Sciences\\University of
Arkansas\\\newline Fayetteville\\AR 72701\\USA}
\email{yoav@uark.edu}
\urladdr{}

\volumenumber{12}
\issuenumber{}
\publicationyear{2007}
\papernumber{10}
\startpage{277}
\endpage{284}

\doi{}
\MR{}
\Zbl{}

\arxivreference{math.GT/0607332}

\received{19 October 2005}
\revised{16 July 2006}
\accepted{16 July 2006}
\published{3 December 2007}
\publishedonline{3 December 2007}
\proposed{}
\seconded{}
\corresponding{}
\version{}

\makeatletter
\def\cnewtheorem#1[#2]#3{\newtheorem{#1}{#3}[section]
\expandafter\let\csname c@#1\endcsname\c@thm}

\makeatother  %  move after \newtheorem block

\let\xysavmatrix\xymatrix
\def\xymatrix{\disablesubscriptcorrection\xysavmatrix}
\AtBeginDocument{\let\bar\wbar\let\tilde\wtilde\let\hat\what}

\newtheorem{thm}{Theorem}
\newtheorem{cor}{Corollary}[section]
\newtheorem{lem}[cor]{Lemma}
\newtheorem{pro}[cor]{Proposition}
\newtheorem{prodef}[cor]{Proposition/Definition}
\newtheorem{lemdef}[cor]{Lemma/Definition}
\theoremstyle{definition}

\theoremstyle{remark}
\newtheorem{rem}[cor]{Remark}
\numberwithin{equation}{section}
\def\s2s{S^2_s}
\def\sst{\Sigma_t}

\def\eee{Esmeralda}
\def\iii{Ivan}
\begin{document}

\begin{asciiabstract}
In [Topology 35 (1996) 1005--1023] J H Rubinstein and M Scharlemann, 
using Cerf Theory, developed tools for comparing Heegaard splittings of 
irreducible, non-Haken manifolds. As a corollary of their work they obtained a 
new proof of Waldhausen's uniqueness of Heegaard splittings of 
S^3. In this note we use Cerf Theory and develop the tools needed 
for comparing Heegaard splittings of S^3. This allows us to use 
Rubinstein and Scharlemann's philosophy and obtain a simpler proof 
of Waldhausen's Theorem. The combinatorics we use are very similar 
to the game Hex and requires that Hex has a winner. The paper 
includes a proof of that fact (Proposition 3.6).
\end{asciiabstract}

\begin{webabstract}
In [Topology 35 (1996) 1005--1023]
J\,H Rubinstein and M Scharlemann, using Cerf Theory, developed 
tools for comparing Heegaard splittings of irreducible, non-Haken 
manifolds. As a corollary of their work they obtained a new proof of 
Waldhausen's uniqueness of Heegaard splittings of $S^3$. In this note 
we use Cerf Theory and develop the tools needed for comparing Heegaard 
splittings of $S^3$. This allows us to use Rubinstein and Scharlemann's 
philosophy and obtain a simpler proof of Waldhausen's Theorem. The 
combinatorics we use are very similar to the game Hex and requires 
that Hex has a winner. The paper includes a proof of that fact 
(Proposition 3.6).
\end{webabstract}

\begin{htmlabstract}
In [Topology 35 (1996) 1005&ndash;1023]
J&nbsp;H Rubinstein and M Scharlemann, using Cerf Theory, developed
tools for comparing Heegaard splittings of irreducible, non-Haken
manifolds. As a corollary of their work they obtained a new proof of
Waldhausen's uniqueness of Heegaard splittings of S<sup>3</sup>. In this note
we use Cerf Theory and develop the tools needed for comparing Heegaard
splittings of S<sup>3</sup>. This allows us to use Rubinstein and Scharlemann's
philosophy and obtain a simpler proof of Waldhausen's Theorem. The
combinatorics we use are very similar to the game Hex and requires
that Hex has a winner. The paper includes a proof of that fact
(Proposition 3.6).
\end{htmlabstract}

\begin{abstract}

In \cite{scharlemann-rubinstein} 
J\,H Rubinstein and M Scharlemann, using Cerf Theory, developed 
tools for comparing Heegaard splittings of irreducible, non-Haken 
manifolds. As a corollary of their work they obtained a new proof of 
Waldhausen's uniqueness of Heegaard splittings of $S^3$. In this note 
we use Cerf Theory and develop the tools needed for comparing Heegaard 
splittings of $S^3$. This allows us to use Rubinstein and Scharlemann's 
philosophy and obtain a simpler proof of Waldhausen's Theorem. The 
combinatorics we use are very similar to the game Hex and requires 
that Hex has a winner. The paper includes a proof of that fact 
(\fullref{pro:esmeralda wins}).
\end{abstract}

\maketitle

\section{Introduction}

In \cite{waldhausen} F Waldhausen proved uniqueness of
non-stabilized Heegaard splittings of $S^3$.

\begin{thm}[Waldhausen]
\label{thm:waldhausen}%
Let $\Sigma \subset S^3$ be a Heegaard surface of genus $g > 0$.
Then $\Sigma$ is a stabilization of a Heegaard surface of genus
$g-1$.
\end{thm}

In \cite{scharlemann-rubinstein} J~H Rubinstein and M Scharlemann,
using Cerf Theory \cite{cerf}, developed tools for comparing Heegaard
splittings of irreducible,
non-Haken manifolds.  As  a corollary of their work they obtained a new
proof
of \fullref{thm:waldhausen}.  In this note we use Cerf Theory and
develop the tools needed for comparing Heegaard splittings of $S^3$.  
This allows us
to use Rubinstein and Scharlemann's philosophy and obtain a simpler proof
of
\fullref{thm:waldhausen}.  We assume familiarity with the basic facts
and
standard terminology of 3--manifold topology and in particular Heegaard
splittings; see Scharlemann \cite{scharlemann-survey}.  For another proof of
Waldhausen's
Theorem see Johnson \cite{johnson}.

We begin with an outline of the proof.  As with many proofs of
\fullref{thm:waldhausen} we assume the theorem is false and pick
$\Sigma$
to be a minimal genus counterexample; we induct on $g$, the genus of
$\Sigma$.
A simple application of van Kampen's theorem shows that if $g=1$ then the
meridians of the complementary solid tori intersect minimally once and
hence
$\Sigma$ is a stabilization of the genus zero splitting of $S^3$. The heart
of
the argument (in the following sections) is to show that if $g > 1$ then
$\Sigma$ weakly reduces.  By A Casson and C McA  Gordon's seminal work
\cite{casson-gordon} either $\Sigma$ reduces or $S^3$ contains an essential
surface.  As the latter is impossible, $\Sigma$ must reduce. Cutting $S^3$
open along the reducing sphere we obtain 2 balls (say $B_1$ and $B_2$,
resp.)
and a once punctured surface in each (say $S_1$ and $S_2$, resp.).  We
attach
3--balls to $B_1$ and $B_2$ and cap off $S_1$ and $S_2$ with disks. It is
easy
to see that we obtain two Heegaard splittings of $S^3$, each of positive
genus less than $g$.  By our inductive hypothesis each of them is
stabilized. Hence, $\Sigma$ is stabilized as well.

The remainder of this paper is devoted to showing that if $\Sigma$ is a
Heegaard splitting of $S^3$ of genus $g>1$ then $\Sigma$ weakly reduces.

\subsection*{Acknowledgment} 

I thank the anonymous referee for helpful suggestions.

\section{The Graphic}

\label{sec:graphic}

$S^3$ is the unit sphere in $\mathbb R^4$.  As such, it inherits a height
function given by the projection onto the $x$--axis, denoted $h_1$.  $S^3$
has
one maximum at $(0,0,0,1)$, one minimum at $(0,0,0,-1)$ and for any $s \in
(-1,1)$ we have that $h_1^{-1}(s)$ is a 2--sphere which we denote $\s2s$.
This
is a special case of a {\em sweepout\/}.

Given $\Sigma$, we have a sweepout of $S^3$ corresponding to $\Sigma$; this
concept was originally introduced in 
Scharlemann--Rubinstein
\cite{scharlemann-rubinstein}.
Although
our description is a little different from that given in
\cite{scharlemann-rubinstein} it is easy to see that the two are
equivalent;
for a more detailed treatment similar to this paper, see
Rieck \cite{rieck-topology,rieck-kjm} and Rieck--Rubinstein
\cite{rieck-rubinstein}.   Let $S^3 =
U
\cup_\Sigma V$ be the Heegaard splitting corresponding to $\Sigma$.  Let
$h_2$
be a height function on $U$, $h_2 \co U \to [-1,0]$ so that $\partial U =
\Sigma$
is at level $0$, a spine of $U$ is at the level -1 and for each $t \in
(-1,0]$, $h_2^{-1}(t)$ is a surface parallel to $\partial U$.  Similarly
take
a height function on $V$ (also denoted $h_2$) $h_2 \co V \to [0,1]$, so that
$\partial V = \Sigma$ is at level $0$, a spine of $V$ is at the level 1 and
for each $t \in [0,1)$, $h_2^{-1}(t)$ is a surface parallel to $\partial
V$. Pasting the two functions together and obtain a function 
$h_2 \co S^3 \to [-1,1]$. For $t \in (-1,1)$ we denote $h_2^{-1}(t)$ by 
$\sst$. By
transversality we may assume that the spines of $U$ and $V$  are disjoint
from
$(0,0,0,1)$ and $(0,0,0,-1)$.  For every point $(s,t) \in (-1,1) \times
(-1,1)$ we have the two surfaces $\s2s$ and $\sst$. Cerf Theory says that
we
can perturb $h_1$ and $h_2$ so that the intersection of $\s2s$ and $\sst$
is
transverse for almost all $(s,t) \in [-1,1] \times [-1,1]$, and the set for
which the intersection is not transverse forms a finite graph (called the
{\em
Graphic\/}) with the following properties.
\begin{enumerate}
\item For $(s,t)$ on an edge of the graphic, $\s2s \cap \sst$
contains exactly one non-degenerate critical point (either
center or a saddle).
\item At a valence 4 vertex the corresponding surfaces have
exactly two non-degenerate critical points.  A valence 4 vertex
can be seen as a point where two arcs of the graphic cross
each other, each corresponding to a single
non-degenerate critical point.
\item There is one other type of vertex (called a \em
Birth-Death\em\ vertex) that has valence 2.  Birth-death
vertices do not play a role in our study and we will not
describe them here.
\end{enumerate}

The closure of a component of $[-1,1] \times [-1,1]$
cut open along the Graphic is called a {\em region\/}.  Given a
region, the intersection of the surfaces that correspond to a
point in the region does not depend in the choice of point in any
essential way.

\section[The labels I, E and a friendly game of Hex]
{The labels $I$, $E$ and a friendly game of Hex}
\label{sec:hax}

We label the regions.  A region is labeled $E$ (standing for ``essential'')
whenever the intersection of surfaces corresponding to a
point in the region contains a curve that is essential in $\sst$;
otherwise, the label $I$ (standing for ``inessential'') is used.  By
definition each region has exactly one label.

\begin{figure}
\cl{\includegraphics[width=2in]{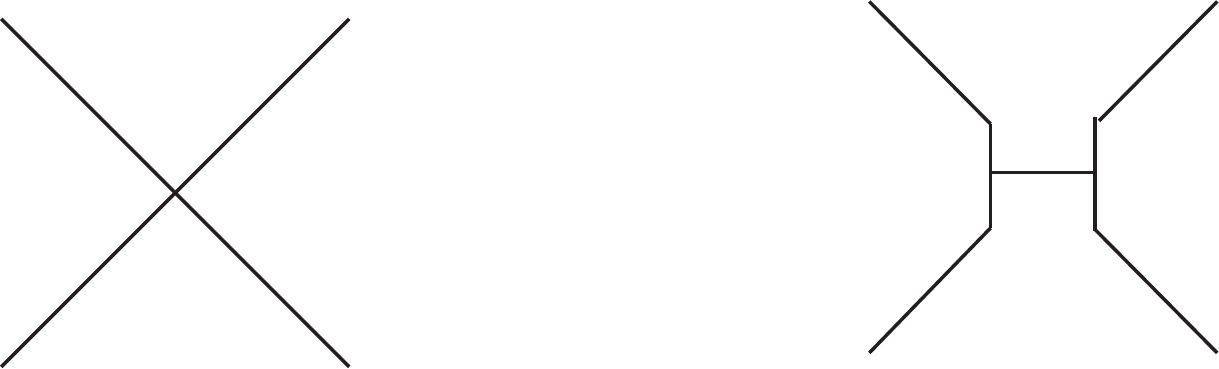}}
\caption{Obtaining the board}
\label{fig:board}
\end{figure}

In order to enjoy a game of Hex we modify the Graphic as follow:
if a valence 4 vertex is adjacent to two $E$--regions and two
$I$--regions and the labels alternate when going cyclically around
it, we split the Graphic and introduce a short edge separating the
$I$ regions; see \fullref{fig:board} where the northern and
southern regions are $I$--regions and the western and eastern regions are
$E$--regions. The graph obtained is called the {\em Board\/}. Note
that there is a natural correspondence between regions of the
Graphic and those of the Board; using this correspondence the
regions of the Board inherit labels from the Graphic.

The reason for creating the Board is the following.

\begin{prodef} \label{pro:boundary is a 2-manifold}%

By the {\em border\/} we mean the union of the edges of the Board that
separate $E$--regions from $I$--regions. In $(-1,1) \times (-1,1)$ the border
forms an embedded 1--manifold.
\end{prodef}

\begin{proof}
Away from the vertices the proposition clearly holds. Let $v$ be a valence
4
vertex.  If all regions around $v$ have the same label $v$ isn't on the
border.  If one region has one label and three have the other label, the
border is locally an interval (with a corner).  If two regions adjacent to
$v$
are labeled $E$ and two are labeled $I$, then by construction of the Board
the
labels do not alternate. Hence the Border cuts across such a vertex
smoothly. At a valence 3 vertex $v$ either all three regions have the same
label (and $v$ is not on the Border) or two regions with one label meet a
third region with the other label (and the border is locally an
interval).
\end{proof}

\begin{rem}
\label{rem:hex}
In the original game of Hex every vertex has valence 3.  Hence the border
there forms a 1--manifold as well.
\end{rem}

We now pick two volunteers to play Hex.  The goal of the first volunteer,
\iii, is  finding a chain of regions (say $R_1,\dots,R_n,\, n \geq 1$)
labeled
$I$ that connects the left edge of the Board (points with $s=-1$) with its
right ($s=1$).  Similarly, the goal of the second volunteer, \eee, is
finding
a chain of regions labeled $E$ that connects the bottom edge of the board
($t=-1$) with its top ($t=1$).  In both cases, the region $R_{i-1}$  shares
an
edge with $R_i$ ($i=2,\dots,n$).

The next proposition is quite special to $S^3$.

\begin{pro}
\label{pro:ivan doesn't win}%
\iii\ can't win.
\end{pro}

\begin{proof}
Suppose \iii\ wins and let $R_1,\dots,R_n$ (for some $n$) be a chain of
regions, starting at the left ($s=-1$) and ending at the right ($s=1$)
(note
that it is possible that $R_1$ meets the left edge in a single point only,
and
similarly for $R_n$ and the right edge).  Consider the corresponding
regions in
the Graphic (still denoted $R_i$). The cost: since some edges of the Board
are
crushed, it is now possible that $R_{i-1}$ shares only a valence 4 vertex
with
$R_{i}$. Given $s \in [-1,1]$ we color $h_1^{-1}([-1,s])$ yellow and
$h_1^{-1}([s,1])$ green.

The proof of the following lemma is an easy innermost disk argument and is
left to the reader.

\begin{lemdef}[regarding $I$--regions]
\label{lem:regarding I-regions}%
If $(s,t)$ is in an $I$--region then the entire surface $\sst$ (except
perhaps
for parts contained in a disk) is either yellow or green (resp.); we say
that
$\sst$ is {\em essentially yellow\/} ({\em green\/} resp.).
\end{lemdef}

We replace the labels $I$ by labels $I(G)$ and $I(Y)$ as follows:
$I$--regions
with essentially green surfaces are labeled $I(G)$ and $I$--regions with
essentially yellow surfaces are labeled $I(Y)$. Of course no surface is
essentially green and essentially yellow simultaneously; this, together
with
\fullref{lem:regarding I-regions}, establishes that every $I$--region gets
exactly one label. In addition, it is easy to see that $I$--regions with $s$
very close to -1 are labeled $I(G)$ and $I$--regions with $s$ very close to
1
are labeled $I(Y)$.  Considering the chain of $I$-regions $R_1,\dots,R_n$,
we
see that $R_1$ is labeled $I(G)$ and $R_n$ is labeled $I(Y)$. Let $i$ be
the
first index with $R_i$ labeled $I(Y)$. Thus $R_{i-1}$ is labeled $I(G)$ and
$R_{i}$ is labeled $I(Y)$.  If $R_{i-1}$ and $R_{i}$ share an edge then
passing from one to the other we cross a single critical point, either a
center or a saddle.  In either case, no essential curve is introduced or
removed (recall we are crossing from one $I$--region to another) and
therefore
labels cannot change. Thus we may assume that we cross a valence 4 vertex
(say
$v$), corresponding to 2 singular points (say $s_1$ and $s_2$).  By
construction of the Board $v$ was obtained from pinching an edge of the
Board
and the remaining two regions adjacent to $v$ are $E$--regions (recall
\fullref{fig:board}); since crossing a center doesn't change an
$I$--region
to an $E$--region we see that both $s_1$ and $s_2$ are saddles.

Moving out of $R_{i-1}$ by crossing $s_1$, we arrive at a region labeled
$E$;
thus crossing the saddle has the effect of changing a single inessential
curve
into two parallel essential curves bounding an annulus.  Since the surface
was
essentially green prior to crossing $s_1$, the annulus between the parallel
curves is essentially yellow (ie the annulus is yellow except
perhaps
for regions contained in a disk). Crossing $s_2$ into $R_{i}$ the label
becomes $I$; hence the 2 parallel curves are pinched together to become a
single inessential curve.  If the pinching is done inside the essentially
yellow annulus (thus turning it into a disk) the surface becomes
essentially
green; hence the boundary of the annulus is pinched  outside the
essentially
yellow annulus.  We obtain an essentially yellow once-punctured torus $T$
or
pair of pants $P$.  Since $R_{i}$ is an $I$-region the boundary of $T$
(resp. $P$) is inessential; hence $g=1$ (resp. $g=0$), contradicting our
assumptions. This establishes \fullref{pro:ivan doesn't win}.
\end{proof}

\begin{rem}[About the game Hex]
In the following proposition we prove that Hex has a winner.  The only
properties of the game we are using are: (1) the Border is  an embedded
1--manifold and (2) the four corners of the Board are adjacent to exactly
one
region each. It is easy to see that these conditions hold for the
traditional
game Hex (recall \fullref{rem:hex}), hence the proof of
\fullref{pro:esmeralda wins} shows the well-known fact that in that
game too there is a winner.
\end{rem}

\begin{pro}
\label{pro:esmeralda wins}%
Hex has a winner.
\end{pro}

\begin{proof}
(We work on the Board.)  First observe that there is exactly one region
adjacent to each corner of $[-1,1] \times [-1,1]$ (these regions correspond
to
disjoint surfaces $\s2s$ and $\sst$).

By \fullref{pro:boundary is a 2-manifold} the Border forms an
embedded
1--manifold in $(-1,1) \times (-1,1)$ and therefore has four types of
components (we note that distinct components of the Border may share a
point
on the boundary):
\begin{enumerate}
\item simple closed curves,
\item arcs connecting an edge to itself,
\item arcs connecting an edge to an adjacent edge and
\item arcs connecting an edge to an opposite edge.
\end{enumerate}

We note that curves of type (1) do not play a role in the proof.  Suppose
there is a arc of type (4), say  connecting the top edge to the bottom
edge.
Then on one side of that arc the regions are all labeled $E$ and therefore
\eee\ wins.  Similarly, if there is an arc connecting left edge to the
right
edge \iii\ wins.  We may therefore assume there are no arcs of type (4).
In
that case, by induction on the arcs of type (2) and (3), we can easily
prove
that some region $R$ is adjacent to all four edges.  If $R$ is labeled $I$
then \iii\ wins and if it is labeled $E$ then \eee\ wins.
\end{proof}

By \fullref{pro:ivan doesn't win} and \fullref{pro:esmeralda wins}
\eee\
wins.  \eee's victory is given by a path of regions in the Board, say
$R_1,\dots,R_n,\, n \geq 1$, all labeled $E$ and connecting the bottom of
the
board to its top.  We consider the corresponding regions in the Graphic,
still
denoted $R_1,\dots,R_n$.  Observe that by construction of the Board,
$R_{i-1}$
still shares an edge with $R_{i}$ ($i=2,\dots,n$). (As in
\fullref{pro:ivan doesn't win} $R_1$ ($R_n$ resp.) may have only
one
point on the bottom edge (top edge resp.).)

\section{The weak reduction}

We complete the proof by finding a weak reduction; this is  a standard
argument in Cerf Theory, originally due to Rubinstein and Scharlemann
\cite{scharlemann-rubinstein}. Denote the handlebodies obtained by cutting
$M$
along $\sst$ by $U_t$ and $V_t$.  First we show the following Lemma.

\begin{lem}[regarding $E$--regions]
\label{lem:meridians in R_i's}%
Let $\s2s$ and $\sst$ be surfaces corresponding to a region labeled $E$.
Then
some curve of $\s2s \cap \sst$ bounds a meridian disk in $U_t$ or $V_t$.
\end{lem}

\begin{proof}
Let $\mathcal{C} \subset \s2s \cap \sst$ denote the collection of curves of
$\s2s \cap \sst$ that are essential in $\sst$.   Since the label is $E$,
$\mathcal{C} \neq \emptyset$.  Consider $\mathcal{C}$ as an embedded
1-manifold in $\s2s$ and  let $c \subset \mathcal{C}$ be an innermost curve
and $D \subset \s2s$ the innermost disk it bounds.  Then in its interior
$D$
intersects $\sst$ in a (possibly empty) collection of curves that are
inessential in $\sst$; a standard disk swap argument gets a disk disjoint
from
$\sst$ with boundary $c$.
\end{proof}

Since $R_1$ contains points with $t$ arbitrarily close to -1 (where $\sst$
collapses to a spine of $U_t$) it is easy to see that curves of
$\s2s\cap\sst$
bound meridians of $U_t$; likewise curves of $\s2s\cap\sst$ in $R_n$ bound
meridians of $V_t$. By \fullref{lem:meridians in R_i's} every region
$R_i$
has  curves of $\s2s\cap\sst$ that bound meridians of $U_t$ or $V_t$.  Let
$i$
be the lowest index so that $R_{i}$ has a curve of $\s2s\cap\sst$ bounds a
meridian of $V_t$. We arrive at the following dichotomy.
\begin{enumerate}
\item ($i=1$) Surfaces corresponding to $R_i$ contain a curve of
$\s2s\cap\sst$ that bounds a meridian in $U_t$ and a curve that bounds
a meridian in $V_t$.
\item ($i>1$) Surfaces corresponding to $R_{i-1}$ contain a curve of
$\s2s\cap\sst$ that bounds a meridian in $U_t$ and surfaces
corresponding to $R_i$ contain a curve of $\s2s\cap\sst$ that bounds a
meridian in $V_t$.
\end{enumerate}

In Case (1) we directly see a weak reduction or reduction (if both meridian
disks bound the same curve).

In Case (2), we note that crossing from $R_{i-1}$ to $R_{i}$ corresponds to
crossing one critical point, either a saddle or a center.  In either case
the
set of essential curves in $\s2s\cap\sst$ corresponding to $R_{i-1}$ can be
isotoped to be disjoint from those corresponding to $R_{i}$; hence $\Sigma$
reduces or weakly reduces.

Since reduction implies a weak reduction, we find a weak reduction in every
case above.  This completes the proof of \fullref{thm:waldhausen}.

\begin{rem}
If $g=1$ it is clearly impossible to find a weak reduction.  Reading
through
the proof, we find exactly one place where the assumption $g>1$ was used:
in
the proof that \iii\ can't win (\fullref{pro:ivan doesn't win}). We
conclude that if we run the Cerf-theoretic argument in that case, it is
actually \iii\ who wins and \eee\ who loses.
\end{rem}

As a concluding remark we mention that it is quite possible that Waldhausen
never intended to study Heegaard splittings of $S^3$, but rather prove the
Poincar\'e Conjecture.  If we replace $S^3$ with a homotopy 3--sphere the
argument above fails miserably, since the ``weak reduction'' we will obtain
consists of immersed disks (small problem, in light of Papakyriakopoulos's
work) that might intersect each other (and hence will not give a weak
reduction at all, even if each disk is embedded).  Even if these problems
are
miraculously overcome, the best we can hope for is a ``reduction'' of the
Heegaard surface via an immersed sphere that intersects the Heegaard
surface
in a single (probably not simple) closed curve.  This is equivalent to the
following condition: {\em the intersection of the  kernels of the two maps
induced on the fundamental group of $\Sigma$ by its inclusion into $U_t$
and
$V_t$ is non-trivial\/}.  This is apparently not the right way to go: it
was
proven by J\,R Stallings in his paper ``How not to prove the Poincar\'e
conjecture'' \cite{stalling}.

% ----------------------------------------------------------------

\bibliographystyle{gtart}
\bibliography{link}

\begin{thebibliography}{}
\providecommand\bibmarginpar{\leavevmode\marginpar}
\def\urlstyle#1{{\tt #1}}

\bibitem{casson-gordon}
\textbf{A\,J Casson}, \textbf{C\,M Gordon},
  \href{http://dx.doi.org/10.1016/0166-8641(87)90092-7} {\emph{Reducing
  {H}eegaard splittings}}, Topology Appl. 27 (1987) 275--283 \xox{MR}{918537}

\bibitem{cerf}
\textbf{J Cerf}, \emph{Sur les diff\'eomorphismes de la sph\`ere de dimension
  trois {$(\Gamma \sb{4}=0)$}}, Lecture Notes in Mathematics, No. 53, Springer,
  Berlin (1968) \xox{MR}{0229250}

\bibitem{johnson}
\textbf{J Johnson}, \href{http://dx.doi.org/10.2140/agt.2005.5.1573}
  {\emph{Locally unknotted spines of {H}eegaard splittings}}, Algebr. Geom.
  Topol. 5 (2005) 1573--1584 \xox{MR}{2186110}

\bibitem{rieck-topology}
\textbf{Y Rieck}, \href{http://dx.doi.org/10.1016/S0040-9383(99)00026-9}
  {\emph{Heegaard structures of manifolds in the {D}ehn filling space}},
  Topology 39 (2000) 619--641 \xox{MR}{1746912}

\bibitem{rieck-kjm}
\textbf{Y Rieck}, \emph{An-annular complexes in 3--manifolds}, Kyungpook Math.
  J. 45 (2005) 549--559 \xox{MR}{2205956}

\bibitem{rieck-rubinstein}
\textbf{Y Rieck}, \textbf{J\,H Rubinstein}, \emph{Invariant {H}eegaard surfaces
  in manifolds with involutions and the {H}eegaard genus of double covers}
  (2006) \xox{arXiv}{math.GT/0607145}

\bibitem{scharlemann-rubinstein}
\textbf{H Rubinstein}, \textbf{M Scharlemann},
  \href{http://dx.doi.org/10.1016/0040-9383(95)00055-0} {\emph{Comparing
  {H}eegaard splittings of non-{H}aken {$3$}--manifolds}}, Topology 35 (1996)
  1005--1026 \xox{MR}{1404921}

\bibitem{scharlemann-survey}
\textbf{M Scharlemann}, \emph{Heegaard splittings of compact 3--manifolds},
  from: ``Handbook of geometric topology'', North-Holland, Amsterdam (2002)
  921--953 \xox{MR}{1886684} \xox{arXiv}{math.GT/0007144}

\bibitem{stalling}
\textbf{J\,R Stallings}, \emph{How not to prove the {P}oincar\'e conjecture},
  Ann. of Math. Stud. 60 (1966) 83--88\;available at
  \href{http://math.berkeley.edu/~stall} {{\tt
  http://math.berkeley.edu/\char'176stall}}

\bibitem{waldhausen}
\textbf{F Waldhausen}, \href{http://dx.doi.org/10.1016/0040-9383(68)90027-X}
  {\emph{Heegaard--{Z}erlegungen der {$3$}--{S}ph\"are}}, Topology 7 (1968)
  195--203 \xox{MR}{0227992}

\end{thebibliography}

\end{document}